\documentclass[proc]{edpsmath}
\usepackage[utf8]{inputenc}
\usepackage{amsmath}
\usepackage{amsfonts}
\usepackage{amssymb}
\usepackage{graphicx}
\usepackage{cite}
\usepackage{hyperref}
\newcommand{\Var}{\mathbf{Var}}
\newcommand{\E}{\mathbf{E}}
\newcommand{\xx}{\mathbf{X}}

\newcommand{\R}{\mathbb{R}}
\newcommand{\ud}{\,\text{d}}
\providecommand{\abs}[1]{\left\lvert#1\right\rvert}
\providecommand{\norm}[1]{\left\Vert#1\right\Vert}
\DeclareMathOperator*{\argmin}{argmin}
\begin{document}
\title{Certified metamodels for sensitivity indices estimation}
\author{Alexandre Janon} \address{Laboratoire Jean Kuntzman, Université Joseph Fourier (Grenoble I), INRIA/MOISE, 51 rue des Mathématiques, BP 53, 38041 Grenoble cedex 9, France}
\author{Maëlle Nodet} \sameaddress{1}
\author{Clémentine Prieur} \sameaddress{1}

  \begin{abstract}
Global sensitivity analysis of a numerical code, more specifically 
estimation of  Sobol indices associated with  input variables, generally
requires a large number of model runs. When those demand too much
computation time, it is necessary to use a reduced model (metamodel) to perform sensitivity analysis, whose outputs are numerically close to the ones of
the original model, while being much faster to run. In this case, estimated
indices are subject to two kinds of errors: sampling error, caused by the
computation of the integrals appearing in the definition of the Sobol indices by
a Monte-Carlo method, and metamodel error, caused by the replacement of the
original model by the metamodel. In cases where we have certified bounds for the
metamodel error, we propose a method to quantify both types of error, and we
compute confidence intervals for first-order Sobol indices.
\end{abstract}
\begin{resume}
 L'analyse de sensibilité globale d'un modèle numérique, plus précisément
l'estimation des indices de Sobol associés aux variables d'entrée, nécessite
généralement un nombre important d'exécutions du modèle à analyser. Lorsque
celles-ci requièrent un temps de calcul important, il est judicieux d'effectuer
l'analyse de sensibilité sur un modèle réduit (ou métamodèle), fournissant des
sorties numériquement proches du modèle original mais pour un coût nettement
inférieur. Les indices estimés sont alors entâchés de deux sortes d'erreur:
l'erreur d'échantillonnage, causée par l'estimation des intégrales définissant
les indices de Sobol par une méthode de Monte-Carlo, et l'erreur de métamodèle,
liée au remplacement du modèle original par le métamodèle. Lorsque nous
disposons de bornes d'erreurs certifiées pour le métamodèle, nous proposons une méthode
pour quantifier les deux types d'erreurs et fournir des intervalles de confiance pour les
indices de Sobol du premier ordre.
\end{resume}
\maketitle

\section{Context}
\subsection{Monte-Carlo estimation of first-order Sobol indices}
Let $Y=f(X_1, \ldots, X_p)$ be our (scalar) output of interest, where the input
variables $X_1,\ldots,X_p$ are modelised as independent random variables of
known distribution.  For $i=1,\ldots,p$, we recall the first-order Sobol index:
\[ S_i = \frac{ \Var( \E( Y | X_i) ) }{ \Var (Y) } \] 
which measures, on a scale of $0$ to $1$, the fraction of the total variability 
of the output caused by the variability in $X_i$ alone.

As $f$ is generally implicitly known ($f$ can \emph{e.g.} be a functional of a
solution of a partial differential equation parametrized by functions of the
input variables $X_1,\ldots, X_p$), one has no analytical expression for $S_i$
and has to resort to numerical estimation. The variances in the definition of
$S_i$ can be expressed as multidimensional integrals over the input parameters
space. We use Monte-Carlo estimates for multidimensional integrals: let 
$\{ \xx^k \}_{k=1,\ldots,N}$ and $\{ \xx'^k \}_{k=1,\ldots,N}$ be two
independent, identically distributed vector samples of
$\xx=(X_1,\ldots,X_p)$. For $k=1,\ldots,N$, we note: 
\[ y_k = f(\xx^k)	\;\;\;\text{ and }\;\;\; y_k' = f(X'^k_1, \ldots, X'^k_{i-1}, X^k_i,	X'^k_{i+1}, \ldots, X'^k_p). \] 
We take the following statistical estimator for $S_i$, introduced in \cite{homma1996importance}:  
\[ \widehat S_i(\mathcal E) = \frac{ \frac{1}{N}\sum_{k=1}^N y_k y_k' -
\left(\frac{1}{N}\sum_{k=1}^N y_k \right) \left( \frac{1}{N}\sum_{k=1}^N y_k'	\right) }
{ \frac{1}{N}\sum_{k=1}^N \left(y_k\right)^2 -	\left(\frac{1}{N}\sum_{k=1}^N y_k\right)^2 } \] 
where $\mathcal E=\left( \{\xx^k\}_{k=1,\ldots,N},
\{\xx'^k\}_{k=1,\ldots,N} \right)$ is our couple of samples the estimator depends on.

\subsection{Reduced basis metamodels}
In order to apply the reduced basis metamodelling, we further assume that
the output $f(\xx)$ depends on a function $u(\xx)$ where, for every input $\xx$, 
$u(\xx)$ satisfies a $\xx$-dependent partial differential equation (PDE).

To make things clear, let us consider an example: we take $p=2$, so that
$\xx=(X_1,X_2)$, and take for $u(t,x;X_1,X_2)$ the solution of the following
$(X_1,X_2)$-dependent initial-boundary value problem (viscous Burgers' equation):
\begin{equation}
\label{e:pde}
\left\{ \begin{array}{c}
\frac{\partial u}{\partial t} + \frac{1}{2} \frac{\partial}{\partial x}(u^2) -
\nu \frac{\partial^2 u}{\partial x^2} = 1 \\
u(t=0,x)=u_{0m}^2+5 \sin(0.5x) \;\; \forall x\in[0;1] \\
u(t,x=0)=b_0 \\ u(t,x=1)=b_1
 \end{array} \right.
\end{equation}
where our input parameters are $(X_1,X_2)=(\nu, u_{0m})$, and $b_0$ and $b_1$
are so that we have \emph{compatibility conditions: } $b_0=u_{0m}^2 \text{
and } b_1=u_{0m}^2+5\sin(0.5)$.

This problem can be analyzed by means of the Cole-Hopf substitution (see \cite{hopf1950partial} for instance), which turns the equation of interest into heat equation, leading to an integral representation of $u$ and well-posedness for $u \in  C^0 \left([0,T], H^1(]0,1[) \right)$.

Note that the $x$ symbol denotes the spatial variable $u(\xx)$ depends on, and is
unrelated with the parameters $X_1$ and $X_2$. Our output can be, for instance:
$ f(\xx)=\int_0^T \int_0^1 u(t,x,\xx) \ud x \ud t $.

For a given value of $\xx$, $u(\xx)$ is generally approximated  using a
numerical method, such as the finite-element method. These methods work by
searching for $u(\xx)$ in a linear subspace of high dimension $\mathcal N$; 
this leads to a large linear system (or a succession of linear systems) to solve
for the coefficients of (the approximation of) $u(\xx)$ in a fixed basis of
$X$. This gives what we call the ``full'' discrete solution, that we denote again
by $u(\xx)$. Even if efficient methods have been developed to solve the linear systems
arising from such discretizations, the large number of unkowns that are to be
found is often responsible for large computation times.

The reduced basis method is based on the fact that $\mathcal N$ has to be large
because the basis we expand $u(\xx)$ in does not depend on the PDE problem that
is being solved; hence it is too ``generic'': it can represent well a large
number of functions, but allows much more degrees of freedom than wanted. We split
the computation into two phases: the offline phase, where we seek a ``reduced
space'', whose dimension $n$ is much smaller than $\mathcal N$, and
which is suitable for effectively representing $u(\xx)$ for various values
of the input parameter $\xx$; and the online phase, where, for each required
value of the input parameters, we solve the ``projected'' PDE on the reduced
space. %This time, we solve a linear system (or a succession of linear systems) with
%only $n$ equations instead of $\mathcal N$, to find our reduced solution
%$\widetilde u(\xx)$. 

This method is interesting if we are to solve the PDE for a number of
values of the parameter sufficiently large so that the fixed cost of the offline
phase is cancelled by the gain in marginal cost offered by the online phase
\emph{vs.} the standard discretization. This is often the case with Monte-Carlo
estimations.

One crucial feature of the reduced basis method, which we will rely on
later, is that it provides a \emph{certified error bound}
$\epsilon_u(\xx,t)$, which satisfies ($\norm{\cdot}$ being the usual norm on $L^2([0,1])$):
\[ \norm{ u(\xx;t) - \widetilde u(\xx;t) } \leq \epsilon_u(\xx,t) \;\;\; \forall \xx, \forall t \in[0;T] \]
and, of course, $\epsilon_u$ can be \emph{fully} and \emph{quickly} computed
(with a computation time of the same order of magnitude than the one for
$\widetilde u$). This error bound on $u$ can lead to an error bound $\epsilon$ on the output:
\[ \abs{ f(\xx) - \widetilde f(\xx) } \leq \epsilon(\xx) \;\;\; \forall \xx \]
where $\widetilde f(\xx)$ denotes a functional of the reduced solution.

%Note that the reduced basis method requires additional hypothesis on the
%underlying equation, as well as some handwork to develop the offline/online
%procedures, and the error estimation, for each type of equation and each type
%of output that are to be considered.

One can turn to \cite{nguyen2005certified} for a detailed introduction to the reduced
basis method, and to \cite{janon2011certified} for the extension to the viscous
Burgers equation \eqref{e:pde}.

\section{Construction of combined confidence intervals}

\subsection{Metamodel error}
For a couple of samples $\mathcal E=\left( \{\xx^k\}_{k=1,\ldots,N},
\{\xx'^k\}_{k=1,\ldots,N} \right)$, we can use our metamodel output $\widetilde f$
and our metamodel error bound $\epsilon$ to compute, for $k=1,\ldots,N$:
\[ \widetilde y_k=\widetilde f(\xx^k),\;\; \widetilde y_k'=\widetilde
f(X'^k_1,\ldots,X'^k_{i-1},X^k_i,X'^k_{i+1},\ldots,X'^k_p) \]
and:
\[ \epsilon_k=\epsilon(\xx^k),\;\; \epsilon_k'=\epsilon(X'^k_1,\ldots,X'^k_{i-1},X^k_i,X'^k_{i+1},\ldots,X'^k_p ) \]
In \cite{janon2011confidence}, we show that we can compute rigorous and accurate bounds $\widehat S_i^m$ and $\widehat S_i^M$ depending only on the
$\widetilde y_k, \widetilde y'_k, \epsilon_k$ and $\epsilon'_k$ so that:
\[ \widehat S_i^m(\mathcal E) \leq \widehat S_i(\mathcal E) \leq \widehat
S_i^M(\mathcal E) \]
where $\widehat S_i(\mathcal E)$ is the (unknown) value of the estimator of $S_i$ computed
on the couple of samples $\mathcal E$. We emphasize that, in our approach, the $y_k$ and $y_k'$ are not observed, as no
evaluation of the full model is performed. 

\subsection{Combined confidence intervals}
To take sampling error in account, we use a bootstrap procedure (see \cite{archer1997sensitivity}) on the two bounding estimators $\widehat S_i^m$ and $\widehat S_i^M$. More specifically, we draw $N$ numbers with repetition from $\{1,2,\ldots,N\}$, so as to get a random list $L$. We then get two \emph{bootstrap replications} by computing $\widehat S_i^m$ and $\widehat S_i^M$ using the samples couple $\left( \{\xx^k\}_{k\in L},\{\xx'^k\}_{k\in L} \right)$ instead of $\left( \{\xx^k\}_{k=1,\ldots,N},\{\xx'^k\}_{k=1,\ldots,N} \right)$. We repeat those computations for a fixed number $B$ of times, so as to obtain $B$ couples of replications $S_i^{m,1},\ldots,S_i^{m,B}$ and  $S_i^{M,1},\ldots,S_i^{M,B}$. Now, for a fixed risk level $\alpha\in]0;1[$, let $S_i^{inf}$ and $S_i^{sup}$ be, respectively, the $\alpha/2$ quantile of $S_i^{m,1},\ldots,S_i^{m,B}$ and $S_i^{sup}$ be the $1-\alpha/2$ quantile of $S_i^{M,1},\ldots,S_i^{M,B}$. We take $[S_i^{inf};S_i^{sup}]$ as our combined confidence interval for $S_i$. This confidence interval accounts for both metamodel and sampling error.

\subsection{Choice of sample size and reduced basis size}
\label{choicesample}
Increasing $N$ and/or $n$
will increase the overall time for computation (because of a larger number of
surrogate simulations to perform if $N$ is increased, or, if $n$ is increased,
each surrogate simulation taking more time to complete due to a larger linear
system to solve). However, increase in these parameters will also improve the
precision of the calculation (thanks to reduction in sampling error for
increased $N$, or reduction in metamodel error for increased $n$). In practice,
one wants to estimate sensitivity indices with a given precision (\emph{ie.} to
produce $(1-\alpha)$-level confidence intervals with prescribed length), and has
no \emph{a priori} indication on how to choose $N$ and $n$ to do so. Moreover,
for one given precision, there may be multiple choices of suitable couples
	$(N,n)$, balancing between sampling and metamodel error. We wish to choose
	the best, that is, the one who gives the smallest computation time. 

On the one hand, we evaluate computation time: an analysis of the reduced basis method shows that the most costly operation made during a call to the metamodel is the resolution of a linear system of $n$ equations; this resolution can be done (e.g., by using Gauss' algorithm) with $O(n^3)$ operations. This has to be multiplied by the required number of online evaluations, i.e. the sample size $N$. Hence, we may assume that computation time is proportional to $N \times n^3$.

On the other hand, the mean length of the
$(1-\alpha)$-level confidence intervals for $S_1,\ldots,S_p$ can be written as the
sum of two terms. The first, depending on $N$, accounts for sampling error and
can be modelled as
$ \frac{2 q_\alpha \sigma}{\sqrt{N}}, $
where $\frac{\sigma}{\sqrt N}$ is the standard deviation of $\widehat S_i$ and $q_\alpha$ is an appropriate $\alpha$-dependent quantile of the standard gaussian distribution. The assumption
of $1/\sqrt N$ decay is heuristically deduced from central limit theorem.

The second term, which accounts for metamodel error, is assumed to be of exponential decay when $n$ increases: $C/a^n$, where $C>0$ and $a>1$ are constants. This assumption is backed up by numerical experiments as well as theoretical works \cite{buffa2009apriori}.

We now wish to minimize computation time while keeping a fixed precision $p$:
\begin{equation}\label{e:coptprob} \text{Find } (N^*,n^*) = \argmin_{(N,n)\in\R^+\times\R^+} n^3 \times N \text{ so that }
\frac{2 q_\alpha \sigma}{\sqrt{N}} + \frac{C}{a^n} = p. \end{equation}

The resolution of this problem is an elementary calculus argument. The solution involve the parameters $C$, $a$ and $\sigma$, which can be fitted against confidence interval lengths found during a ``benchmark run''.

\section{Numerical results}
\subsection{Target model}
Our underlying model is given by the Burgers equation \eqref{e:pde}. Our output functional is taken to be:
$ f(\xx)=\frac{1}{\mathcal N} \sum_{i=0}^{\mathcal N} u \left(t=T, x=\frac{i}{\mathcal N} ; \xx \right) $.

We set $\mathcal N=60$, $\Delta t=.01$, $T=.05$, while the uncertain parameters $\nu$ and $u_{0m}$ are assumed to be of uniform distributions, with respective ranges $[1;20]$ and $[-0.3;0.3]$. We also take $B=300$ bootstrap replications and a risk level $\alpha=0.05$.

Note that more flexible parametrizations of right-hand sides in \eqref{e:pde} can be considered; results remain qualitatively the same. We chose this parametrization for simplicity reasons.

\subsection{Convergence benchmark}
\begin{figure}
\begin{center}
\includegraphics[height=6cm]{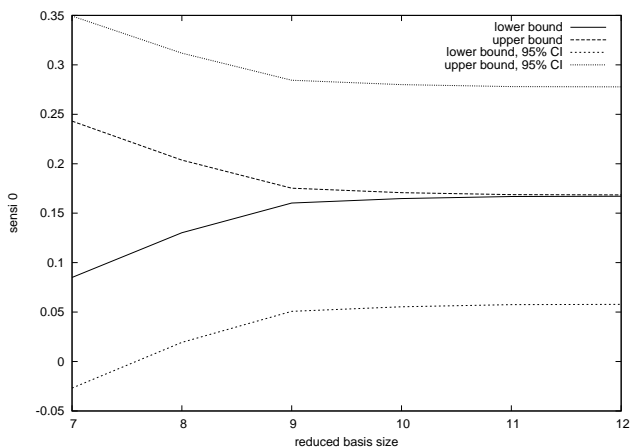} \\
\caption{Convergence benchmark for sensitivity index of $\nu$. We plotted, for a fixed sample of size $N=300$, estimator bounds $\widehat S^m$ and $\widehat S^M$, and endpoints of the 95\% combined confidence interval, for different reduced basis sizes.}
\label{f:fig1}
\end{center}
\end{figure}

Figure \ref{f:fig1} shows the lower $\widehat{S^{m}}$ and upper
$\widehat{S^{M}}$ bounds for different reduced basis sizes $n$ and fixed sample of size $N=300$, as well as the endpoints of the combined confidence intervals. This figure exhibits the fast convergence of our bounds to the true value of $S_a$ as the reduced basis size increases. We also see that the part of the error due to sampling (gaps between confidence interval upper bound and upper bound, and between confidence interval lower bound and lower bound) remains constant, as sample size stays the same.

\subsection{Comparison with estimation on the full model}
To demonstrate the interest of using sensitivity analysis on the reduced model, we computed the combined intervals for the two sensitivity indices using sample size $N=22000$ and $n=11$ (those parameters are found using the procedure described in Section \ref{choicesample} for a target precision $p=0.02$). We found $[0.0674128 ; 0.0939712]$ for sensitivity index for $\nu$, and $[0.914772 ; 0.926563]$ for sensitivity with respect to $u_{0m}$. These confidence intervals have mean length:
$0.019 \approx 0.02 $ as desired. This computation took 58.77 s of CPU time to complete (less than 1 s being spent in the offline phase).

To obtain a result of the same precision, we carry a simulation on the \emph{full} model, for $N=22000$ (sample size can be chosen smaller than before, as there will be no metamodel error); we get a bootstrap confidence interval with mean length of $\approx0.0193$ (we can only provide a confidence interval, as the exact values of the sensitivity indices are not known in this case). This computation takes 294 s of CPU time to complete. Hence, on this example, using a reduced-basis model roughly divides overall computation time by a factor of 5, without any sacrifice on the precision and the rigorousness (as our metamodel error quantification procedure is fully proven and certified) of the confidence interval. We expect higher time savings with more complex (for example, two- or three-dimensional in space) models.  

\begin{acknowledgement}
This work has been partially supported by the French National
Research Agency (ANR) through COSINUS program (project COSTA-BRAVA
no ANR-09-COSI-015).
\end{acknowledgement}

\bibliographystyle{plain}
\bibliography{biblio}

\end{document}